\documentclass{amsart}
\usepackage{epsf} 
\usepackage{amssymb,latexsym, amsmath, amscd, array} 
 
\swapnumbers 

\theoremstyle{plain} 
\newtheorem{thm}{Theorem}
 
\newtheorem{cor}[thm]{Corollary} 
 
\newtheorem{theorem}[thm]{Theorem}

\newcommand\theoref{Theorem~\ref}

\theoremstyle{definition}

\newtheorem{rem}[thm]{Remark}

\def\p{{\noindent \it Proof. }}

\def\ga{\alpha} 
\def\gb{\beta}

\def\R{{\mathbb R}} 
\def\N{{\mathbb N}}

\def\wh{\widehat} 
\def\m{\medskip}

\begin{document} 

\title{Fundamental groups of compact Hausdorff spaces}

\author[J.~Keesling]{James E. Keesling}
\address{J. Keesling, Department of Mathematics, University 
of Florida, 358 
Little Hall, Gainesville, FL 32611-8105, USA} 
\email{jek@math.ufl.edu} 

\author[Yu.~Rudyak]
{Yuli B. Rudyak$^{*}$}\address{Yu. Rudyak, Department of Mathematics, University 
of Florida, 358 
Little Hall, Gainesville, FL 32611-8105, USA} 
\email{rudyak@math.ufl.edu} 
\thanks{$^{*}$Supported by NSF, grant 0406311}

\begin{abstract} We discuss which groups can be realized 
as the fundamental groups of compact Hausdorff spaces. In particular, we prove 
that the claim ``every group can be realized as the fundamental group of a 
compact Hausdorff space'' is consistent with the Zermelo--Fraenkel--Choice set 
theory.
\end{abstract}

\maketitle 

In this note we prove that every group $\pi$ can be realized as the fundamental 
group of a compact Hausdorff space, if we work in the pattern of the set theory 
without inaccessible cardinals. (See e.g. \cite{Kan} concerning inaccessible 
cardinals.) Notice that, because of a theorem of Kuratowski \cite[Prop. 1.2]{Kan}, the absence of inaccessible cardinals is consistent with the Zermelo -- Fraenkel -- Choise set theory.

The idea is the following. 
Take a $CW$-space $X$ with $\pi_1(X)=\pi$ and let $\gb X$ be the Stone--\v Cech 
compactification of $X$, \cite{M}. Then $X$ is a path-connected component of 
$\gb X$, and therefore $\pi_1(\gb X, *)=\pi$ for all $*\in X$.

This result should be contrasted with the result of Saharon Shelah \cite{S} that 
for path connected, locally path connected compact metric spaces $X$, 
$\pi_1(X)$ is either a finitely generated group or has cardinality 
$2^{\aleph_0}$.

\m
We are grateful to Alex Dranishnikov for useful discussions. We must also 
mention that this paper initiated and stimulated by the Algebraic Topology 
Discussion List. 

\m 
All spaces are assumed to be Hausdorff, all maps and functions are assumed to be 
continuous. We denote by $I$ the unit segment $[0,1]$.  Let $\N^*$ denote the 
one-point compactification of the natural numbers $\N$ with $*$ being the point 
at infinity.

A dated, but useful, compendium of information on the Stone--\v Cech 
compactification can be found in \cite{W}.  See also \cite{GJ}.
The  Stone--\v Cech compactification $\gb X$ of the completely regular space $X$ 
can be characterized as a topological embedding $X \subset \gb X$ with $\gb X$ 
compact and such that $X$ is dense in $\gb X$ and such that every function 
$f : X \to I $ can be extended to a function $\wh f: \gb X \to I$. There are 
several standard constructions of the Stone--\v Cech compactification.  One 
construction uses the maximal ideals in the ring of bounded functions on $X$, 
$C^*(X)$, with the hull-kernel topology \cite[Chapter 7]{GJ}.

Let $\upsilon X$ denote the {\it Hewitt realcompactification} of $X$, 
\cite[Chapter  8]{GJ}. Recall that $X\subset v\upsilon X \subset \gb X$ and is 
characterized by the property that every continuous function $X \to \R$ (not 
necessarily bounded) can be extended to $\upsilon X$. A space $X$ is called {\it 
realcompact} if $\upsilon X=X$. Notice that a realcompact space is not 
necessarily compact.

\begin{theorem}\label{realcompact}
 If $X$ is a paracompact space of the non-measurable cardinality then $\upsilon 
X=X$.
\end{theorem}

\p Katetov \cite{K} proved that a paracompact space $X$ is realcompact iff each 
of its closed discrete subspaces is realcompact. On the other hand, a discrete 
subspace is realcompact iff it has non-measurable cardinality, \cite[Chapter 
12]{GJ}
\qed

\begin{theorem}\label{upsilon}
Each nondiscrete, closed subset in $\gb X\setminus \upsilon X$ contains a copy 
of $\gb \N$, and so its cardinality is at least $2^{\mathfrak c}$.
\end{theorem}

\p See \cite[Theorem 9.11]{GJ}.
\qed

\begin{theorem}\label{path}
If $X$ is a path connected  paracompact space of the non-measurable cardinality 
then $X$ is a path component of $\gb X$.
\end{theorem}

\p Suppose that there exists a path $\ga: I \to \gb X$ with $\ga(0)\in X$ and 
$\ga(1)\in \gb X \setminus X$.  Then,  by Theorems \ref{realcompact} and 
\ref{upsilon}, $\ga(I) \cap (\gb X \setminus X)$ is a discrete set.  So, we 
may assume without loss of generality that $\ga([0,1)) \subset X$ and $\ga(1) 
\in \gb(X) \setminus X$. Clearly, then $\ga([0,1))$ is an infinite set: 
otherwise $\ga(1)\in X$.

Now let $\{ t_n \}_{n=1}^\infty$ be a sequence of distinct points in 
$\ga([0,1))$ converging to $\ga(1)$.  Define $f(t_n) = n \mod 2 : \{ t_n 
\}_{n=1}^\infty \to [0,1]$.  Since $X$ is paracompact and therefore normal, 
let $ F : X \to [0,1]$ be any extension of $f$ to all of $X$ using 
Tietze's Extension Theorem.  Then let $\wh 
F : \gb X \to [0,1]$ be the Stone-\v Cech extension of $F$ to $\gb X$.  Then 
$\wh F | \{ t_n \}_{n=1}^\infty = f$ has an extension to $\{ t_n \}_{n=1}^\infty 
\cup \ga(1)$, which is clearly a contradiction.
\qed

\begin{cor}
Every group of the non-measurable cardinality is the fundamental group of a 
compact space.
\end{cor}

\p Let $\pi$ be a given group.  Let $X$ be a connected $CW$-space having 
$\pi_1(X)=\pi$. We can assume that $X$ has the non-measurable cardinality.
Since $X$ is paracompact, we conclude that $\pi_1(\gb X, x_0)=\pi$ 
for any $x_0\in X$ in view of \theoref{path}.
\qed

\begin{rem} Let us recall that every measurable cardinal is inaccessible, 
\cite{Kan}.
\end{rem}

\begin{rem} We still have an open question: Which groups can be realized as the 
fundamental groups of path connected compact Hausdorff spaces.
\end{rem}   

\begin{rem} If we consider measurable cardinals, then \theoref{path} turns out 
to be wrong. However, we still do not have any example of (measurable) group 
that cannot be realized as the fundamental group of a compact space.
\end{rem}


\begin{thebibliography}{99999}

\bibitem{GJ}
L. Gillman and M. Jerison, {\it Rings of Continuous Functions}, Springer-Verlag, 
1976 (Reprint of the 1960 edition.)

\bibitem{Kan}
A. Kanamori, {\it The Higher Infinite}, Springer-Verlag, 1997.

\bibitem{K}
M. Kat\v etov, {\it 
Measures in fully normal spaces}.
Fund. Math. {\bf 38}, (1951) 73--84.

\bibitem{M}
J. Munkres, {\it Topology}, Second Edition; Prentice-Hall, 2000.

\bibitem{S}
S. Shelah, {\it Can the fundamental {\rm (}homotopy{\rm )} group of a space be
the  rationals?} Proc. Amer. Math. Soc. {\bf 103} (1988), 627-632.

\bibitem{W}
R. Walker, {\it The Stone--\v Cech Compactification}, Springer-Verlag, 1974.

\end{thebibliography}
\end{document}